\documentclass[runningheads]{llncs}

\usepackage{epsfig}

\usepackage{amssymb,amsmath} 
\usepackage{latexsym}

\renewcommand{\t}{\theta}
\renewcommand{\a}{\alpha}
\newcommand{\ga}{\gamma}
\newcommand{\T}{{\bf T}}
\newcommand{\RR}{{\mathbb R}}
\begin{document}
\title{Level sets of functions and symmetry sets of smooth surface sections}
\titlerunning{Level sets} 
\author{Andr\'{e} Diatta \inst{1} \and Peter Giblin\inst{1} \and Brendan Guilfoyle\inst{2} \and Wilhelm Klingenberg\inst{3}}
\authorrunning{Diatta,  Giblin,  Guilfoyle and Klingenberg}

\institute{University of Liverpool, adiatta@liv.ac.uk, pjgiblin@liv.ac.uk
 \and Institute of Technology, Tralee, brendan.guilfoyle@ittralee.ie \and University of Durham, 
 wilhelm.klingenberg@durham.ac.uk} 

\maketitle

\begin{abstract} 
We prove that the level sets of a real $C^s$ function of two variables near a non-degenerate critical point are of class 
$C^{[s/2]}$ and apply this to the study of planar sections of surfaces close to the singular section
by the tangent plane at an elliptic or hyperbolic point, and in particular at an umbilic point.

 We also analyse the cases coming from degenerate critical points, corresponding to cusps of Gauss on surfaces, where the differentiability is now reduced to
$C^{[s/4]}$.

 However in all our applications to symmetry sets of planar sections of surfaces, we assume the $C^\infty$ smoothness.
\end{abstract}

\section{Intoduction}

The  medial axis or skeleton of a closed plane curve $\ga$ encode a great deal of information about
the shape of the region enclosed by the curve. This information has been exploited in many ways\footnote{Typing
`medial axis' into an internet search engine produces hundreds of references, to character recognition, Voronoi diagrams, 
interrogation, reconstruction, modification and design of shape, etc.}.  The medial axis can be defined, for
a smooth curve $\ga$, 
as the closure of the locus of centres of `bitangent' circles---tangent to $\ga$ in two places---and whose radius
equals the minimum distance from the centre to $\ga$.  The symmetry set is the closure of the locus of centres
of all bitangent circles. It thus has a richer structure which underlies the transitions of the medial axis.
See for example~\cite{MathsOfSurfaces2000}.

\begin{figure}
  \begin{center}
  \leavevmode
\epsfysize=1.5in
\epsffile{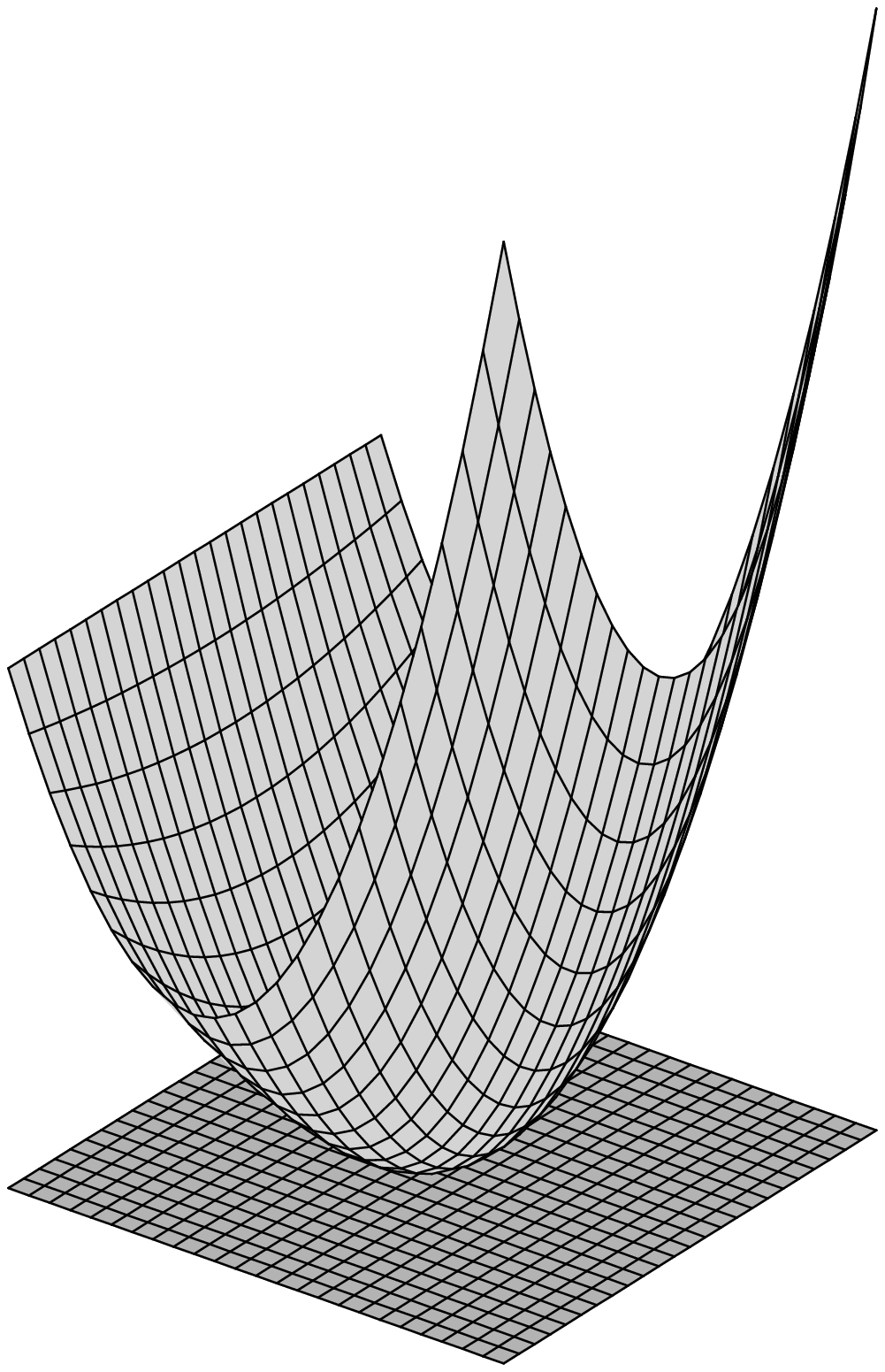}
\epsfysize=1.1in
\epsffile{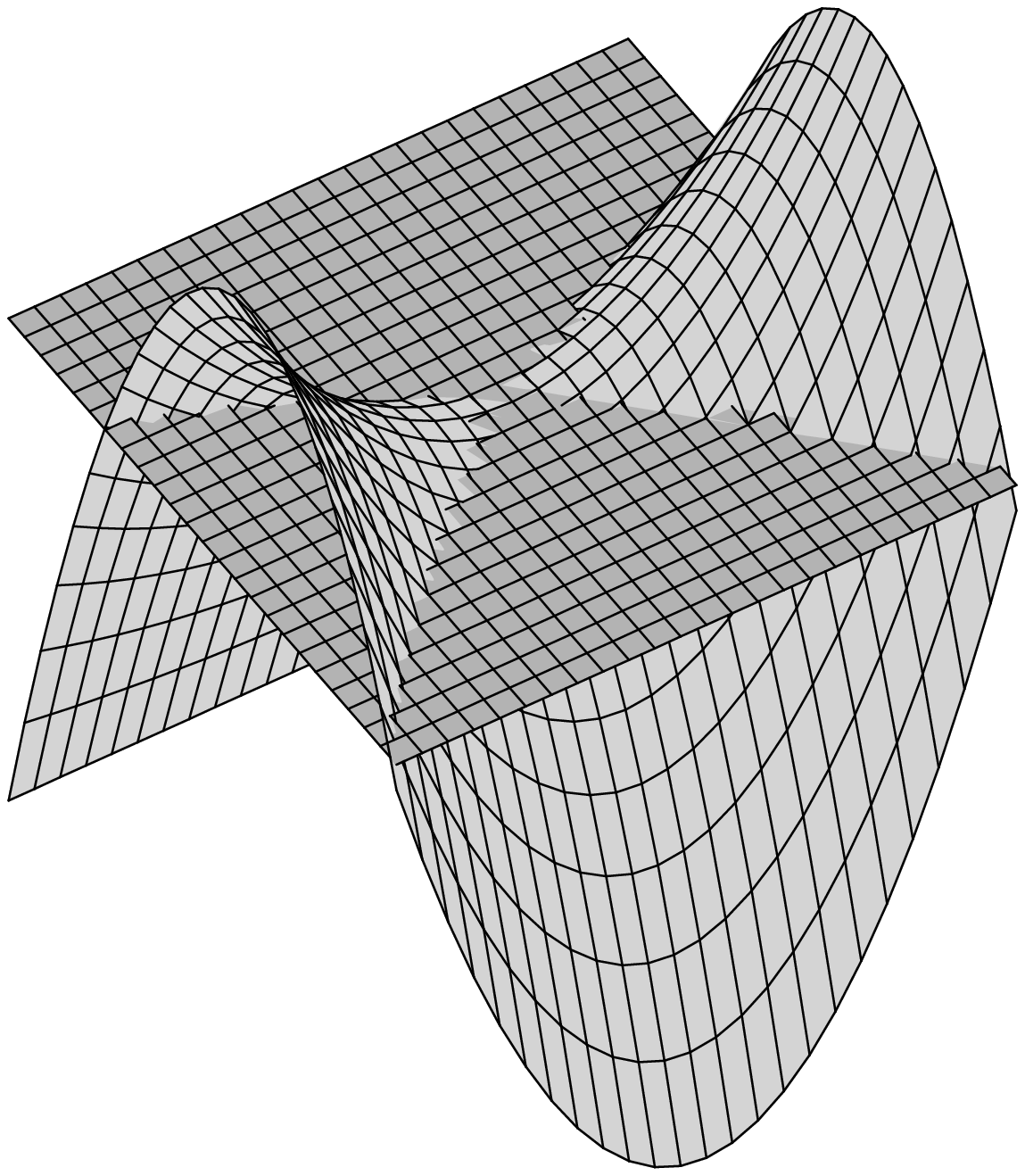}
\epsfysize=1.1in
\epsffile{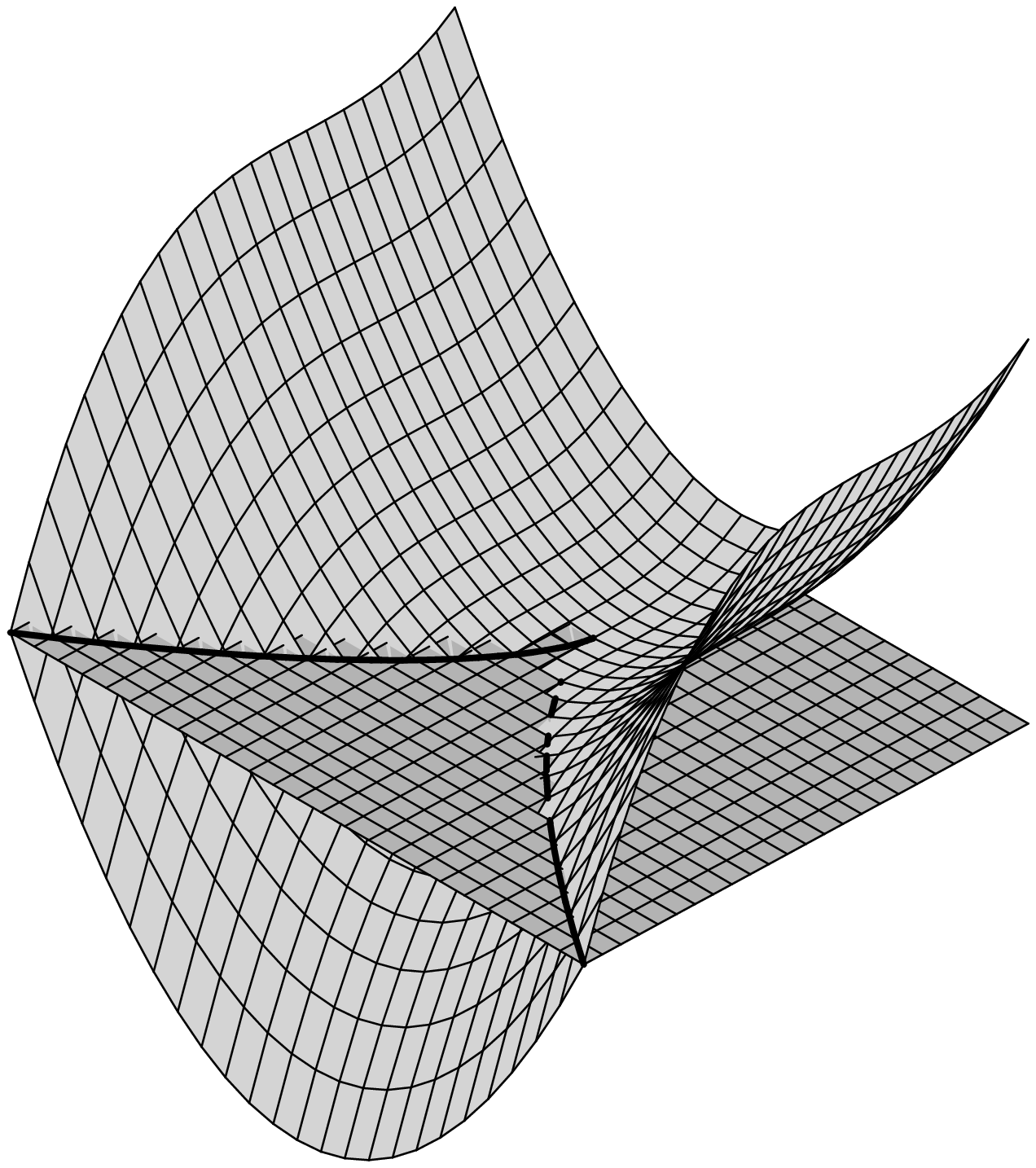}
\epsfysize=1.1in
\epsffile{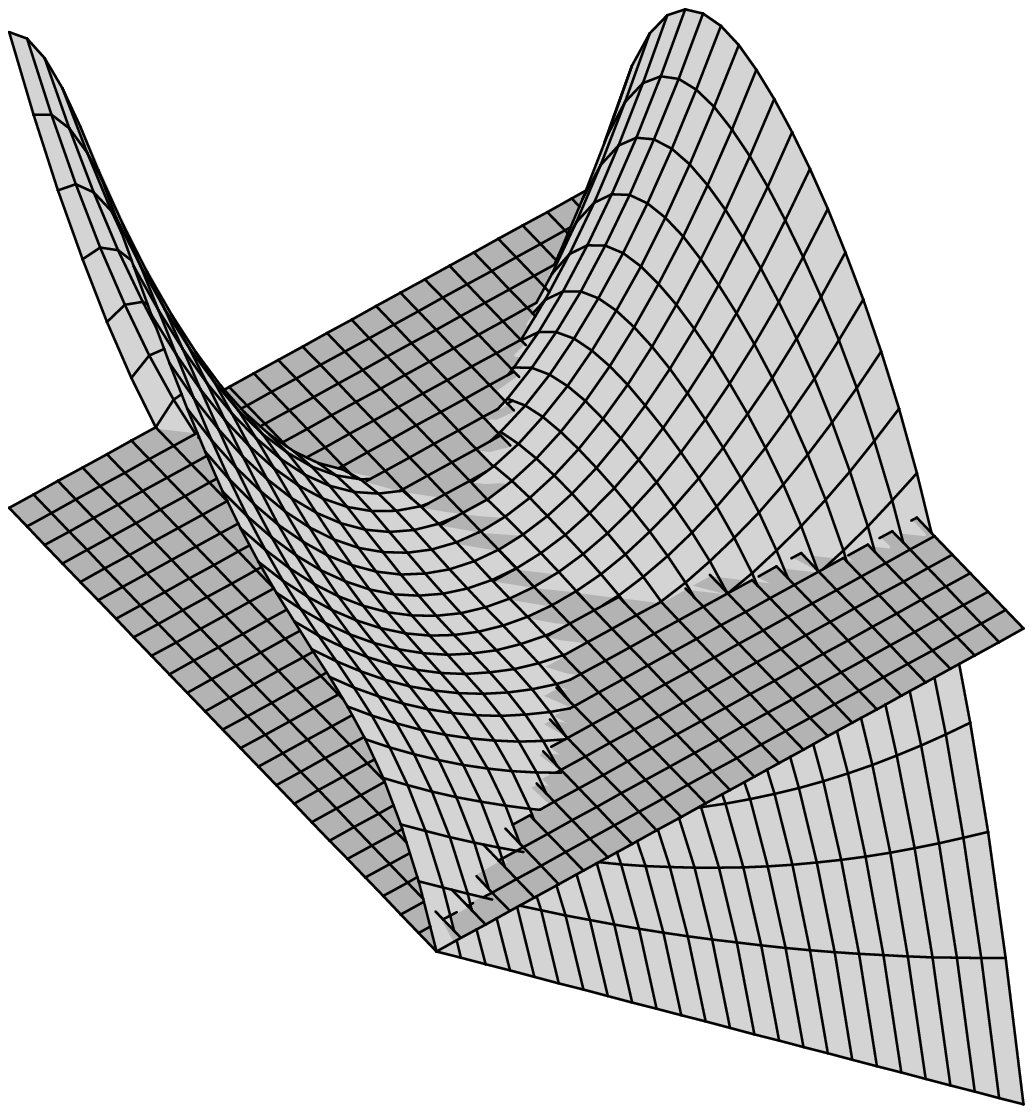}
 \end{center}
\caption{Surfaces with the tangent planes at elliptic, hyperbolic, parabolic and hyperbolic cusp of Gauss points. The intersections
are a point, two smooth transverse curves, a cusped curve and two tangential curves respectively. Figure produced with MAPLE.}
\label{fig:surface-plane}
 \end{figure}

The local structure of the symmetry set of a generic plane curve is of four kinds: smooth branches, endpoints, cusps and triple crossings.
All of these can be seen in the example of Figure~\ref{fig:umbilic}. (For medial axes there are only
smooth branches, endpoints and Y-junctions.)
Symmetry sets and medial axes of curves which vary in a generic 1-parameter family $\ga_k$ are well understood
provided the curves of the family remain nonsingular; a complete local classification was given in \cite{growth-motion}.
In this paper, we consider instead the family of plane sections $\ga_k$ of a smooth surface $M$ in 3-space, where the plane moves
parallel to itself and becomes the tangent plane to $M$ at a point {\bf p} for say $k=0$. Then the section $\ga_0$ 
will be singular at the tangency point {\bf p}. For {\bf p} an elliptic point $\ga_0$ is, locally, a single point.
For {\bf p} a hyperbolic point $\ga_0$ is locally two smooth transversally intersecting branches, while at an ordinary parabolic point
$\ga_0$ has an ordinary cusp at {\bf p}. See Figure~\ref{fig:surface-plane}, which also shows a `hyperbolic cusp of
Gauss' where the intersection is two smooth tangential branches.  See \S\ref{s:cog}. In this situation, the results
of \cite{growth-motion} are invalid and we need to use other techniques to find out how the symmetry sets (and
medial axes) behave in the family of curves.

The motivation for this work comes from the study of isophote curves of a 2D image, which are sections of an intensity surface
in 3-space. Results on the pattern of vertices (extrema of curvature)
and inflexions (zeros of curvature) of the sections $\ga_k$ have been given in \cite{diatta-giblin2004}
while in \cite{scale-space2005} there are results on the patterns of cusps and triple crossings of the symmetry set. Let the
surface $M$ be given locally by an equation $z=f(x,y)$. Then the sections
$\ga_k$ are by their nature given by {\em equations} $f(x,y)=k$, and not by {\em parametrizations} $\ga(t)=(X(t),Y(t))$.
Now plotting symmetry sets of parametrized curves is reasonably straightforward (see \S\ref{s:appl}). But curves given
by equations are a different matter altogether: there is not even any foolproof way to determine how many
real components a curve $f(x,y)=k$ has.  Furthermore it may not be feasible to find an exact parametrization of one
component of a curve $f(x,y)=k$.

In this paper we give a result on parametrizing plane sections of surfaces and use it to approximate such a section to any desired
degree of accuracy. We apply this to the plotting of symmetry sets and medial axes of the surface sections
$\ga_k$. In \S\ref{s:theory} we give the main theoretical result and its proof (Propositions~\ref{prop1}-\ref{prop4}), in \S\ref{s:practice} we give the
way in which this is implemented in practice, and in \S\ref{s:appl} we apply the method to symmetry sets,
concentrating in \S\ref{s:umbilic} on the case of umbilic points (elliptic points where the principal curvatures
are equal). In \S\ref{s:cog} we consider the case of an `elliptic cusp of Gauss', for which the curves $\ga_k$
are simple closed curves when $k>0$.

\smallskip\noindent
{\em Acknowledgements} \ 
The work of the first two named authors 
is a part of the DSSCV project supported by the IST Programme of the European Union (IST-2001-35443).
They also acknowledge EPSRC for the provision of a graphics computer used to draw the symmetry sets.

\section{Level sets of functions}\label{s:theory}
The intuitive idea here that we `blow up' the origin to turn the surface $z=f(x,y)$ into one with nonsingular sections.
We shall take the surface $M$ in Monge form, that is with
\begin{eqnarray}
f(x,y)&=&\textstyle{\frac{1}{2}}\displaystyle (\kappa_1x^2 + \kappa_2y^2)+b_0x^3+b_1x^2y+b_2xy^2+b_3y^3\nonumber\\
&+&c_0x^4+
c_1x^3y+c_2x^2y^2+c_3xy^3+c_4y^4\nonumber \\
&+&d_0x^5+d_1x^4y+
d_2x^3y^2+d_3x^2y^3+d_4xy^4+d_5y^5+ \mbox{ h.o.t}.  
\label{eq:monge}
\end{eqnarray}
For simplicity, consider the case where $\kappa_1$ and $\kappa_2$ are $>0$ so that, locally to the origin,
$z\ge 0$ at points of  $M$. 
Let us write $x=tX, \ y=tY, z=t^2$; after cancellation of $t^2$ the equation $z=f(x,y)$ then becomes
\[1= \textstyle{\frac{1}{2}}\displaystyle (\kappa_1X^2 + \kappa_2Y^2)+t(b_0X^3+\ldots)+t^2(c_0X^4+\ldots)+\ldots,\]
which is a smooth surface in $(X,Y,t)$-space whose sections $t=$ constant $\ne 0$ are scaled versions of the sections
of $M$. We proceed to give the formal details of this procedure in a more general setting and in the next section show how in practice
we have used it.

\medskip

We prove that the level sets of a real function of two variables near a non-degenerate critical point are of class 
C$^{[s/2]-1}$ if the function is of class C$^s$.  To this end we consider $f\in C^s(D)$, $D\subset\RR^2$ the 
open unit disc, $f(0)=df(0)=0$, and set $q=\nabla df(0)$ to be the hessian of $f$ at the origin. We assume that $q$ is either 
positive definite or indefinite.  See Proposition~\ref{prop4} below for the case of a degenerate 
critical point.

\begin{proposition}
Assume that $q$ is positive definite and let $r_0:S^1\rightarrow\RR^+$, where $S^1\subset\RR^2$ is the unit circle, be
defined by $r_0(\theta)=[q(\cos\theta,\sin\theta)]^{-1/2}$, $\theta\in S^1$ corresponding to $e^{i\theta}$. 

Then there exists a function $r:[0,\epsilon)\times S^1\rightarrow\RR^+$ for some $\epsilon>0$ with
\begin{enumerate}
\item[\rm{(a)}] $r\in C^{[s/2]}([0,\epsilon)\times S^1)$,
\item[\rm{(b)}] $r(0,\theta)=r_0(\theta)$ on $S^1$,
\item[\rm{(c)}] $f(t\;r(t,\theta)\cos\theta,t\;r(t,\theta)\sin\theta)=t^2$ on $[0,\epsilon)\times S^1$.
\label{prop:elliptic}
\end{enumerate}
\label{prop1}
\end{proposition}

\noindent{\bf Remark}: The equation (c) implies that for fixed $t\in[0,\epsilon)$, the parametrized curve 
$\theta\mapsto (r(t,\theta)\cos\theta,r(t,\theta)\sin\theta)$ is a rescaled level curve of the
function $f$.

\proof
For $t>0$ consider the curve 
\[
C_t=\{f(tx,ty)=t^2,\; z=t^2\}\subset\RR^3_{xyz}
\]

Consider the surface defined by
\[
Q=\{f(\sqrt{|z|}\;x,\sqrt{|z|}\;y)=|z|\mbox{ for } \ z\neq0 \quad\mbox{and} \quad q(x,y)=1  \mbox{ for }  z=0\}.
\]
We claim that $Q$ is a regular surface of class $C^{[s/2]}$ near $\{z=0\}\cap Q$. Given this, we have 
$C_t=Q\cap\{z=t^2\}$, and since the intersection is transversal (as is seen by computing 
$\mbox{grad }f(\sqrt{|z|}\;x,\sqrt{|z|}\;y)-|z|$ ), we conclude that $C_t\in C^{[s/2]}$ with uniform
bounds. Using the implicit function theorem we derive the existence of $r$ and hence (a) to (c).  
Note that now $C_t$ is parametrized by
\[
\{(t\;r(t,\theta)\cos\theta,t\;r(t,\theta)\sin\theta,t^2) \mbox{ for } \theta\in S^1\}.
\]
 \qed

The claim above follows from the following proposition:

\begin{proposition}
Let $h(x)$ be a real function of one real variable, $h=O(|x|^2)$ and $h\in C^l$. Then $g(x,y)=h(\sqrt{|y|}\;x)$ 
is of class $C^{[l/2]}$ and $g=O(|y|)$.
\label{prop2}
\end{proposition}

\proof
 The proof of this proposition follows from Taylor's formula. \qed

Finally, we have the analogous result for $q$ indefinite, which is proved in a similar way:
\begin{proposition}
Assume that $q$ is indefinite and of the form $q(x,y)=\kappa_1x^2-\kappa_2y^2$ with $\kappa_1\ge\kappa_2>0$.
Let $S^+\subset S^1$ be defined by $q|_{S^+}>0$.
For any $S_0 \subset S^+, \ S_0\ne S^+$ there exists $\epsilon > 0$ and $r :  [0,\epsilon)\times S_0  \to \RR^+$ 
with
\begin{enumerate}
\item[\rm{(a)}] $r\in C^{[s/2]}([0,\epsilon)\times S_0)$,
\item[\rm{(b)}] $f(t\frac{1}{\sqrt{\kappa_1}}r(t,\theta)\cosh\theta,t\frac{1}{\sqrt{\kappa_2}}r(\theta,t)\sinh\theta)=t^2$ 
on $S_0\times[0,\epsilon)$
\item[\rm{(c)}] $r(\theta, 0) = 1$  on $S_0$.    \qed
\end{enumerate}
\label{prop3}
\end{proposition}

We have a result similar to Proposition~\ref{prop1}  in the case of an elliptic cusp of Gauss. In this case
the quadratic terms of $f$ are degenerate, equal to $x^2$ say, so that the surface $z=f(x,y)$ has a
parabolic point at the origin. By making $x$ divide the cubic terms we can still have a closed curve
of intersection $f(x,y)=k>0$ provided the terms in $x^2, xy^2, y^4$ give a positive definite quadratic
form in $x$ and $y^2$. This is called an elliptic cusp of Gauss. The contact of the surface $z=f(x,y)$
with its tangent plane at the origin is of type $A_3$ in the notation of Arnold. See for example
\cite{cog,solid-shape} for many geometrical properties of these points.

\begin{proposition}
Assume that $f \in C^s(D)$ and let the point $0 \in M$ be an elliptic cusp of Gauss
and $\kappa_1 > 0$. Set $q(x,y) = \frac{1}{2}\kappa_1x^2 + b_2xy^2 + c_4y^4$, with
$b_2^2<2\kappa_1c_4$. Let
$r_0 :S^1 \to \RR^+$ be defined by $r_0(\theta) = q(\cos\theta, \sin\theta)^{ -1/4}$.
Then there exists a function $r: [0,\epsilon) \times S^1  \to \RR^+$ for some $\epsilon >0$ with
\begin{enumerate}
\item[\rm{(a)}] $r\in C^{[s/4]}([0,\epsilon)\times S^1)$,
\item[\rm{(b)}] $r(0,\theta)=r_0(\theta)$ on $S^1$,
\item[\rm{(c)}] $f((t \;r(t,\theta))^2\cos\theta,t\;r(t,\theta)\sin\theta)=t^4$ on $[0,\epsilon)\times S^1$.  
\end{enumerate}

\label{prop4}
\end{proposition}

\proof
The rescaling here is nonhomogeneous and given as follows. Let
\[
C_t=\{f(t^2 x,ty)=t^4,\; z=t^4\}\subset\RR^3_{xyz}.
\]
Consider the surface defined by
\[
Q=\{f(|z|^{\frac{1}{2}}\;x,|z|^{\frac{1}{4}}\;y)=|z|\mbox{ for } \ z\neq0 \quad\mbox{and} \quad q(x,y)=1  \mbox{ for }  z=0\}.
\]
Now $Q$ is a regular surface of class $C^{[s/4]}$.This follows from a result analogous to Proposition 2.
\qed

\noindent
{\bf Remark:} \ For a fixed value of $t$, we are parametrizing the level set $f(x,y)=t^4$ not in `polar coordinates'
where each ray from the origin intersects the curve in one point, but by means of parabolas of the form $y^2=kx$ for constants $k$.

\section{Finding the parametrization in practice}\label{s:practice}

All functions and surfaces from now on will be assumed smooth of class $C^\infty$, that is $s=\infty$ in the
Propositions of \S\ref{s:theory}.

We seek to approximate the sections $f(x,y)=k$ of the surface $M$ up to a suitable order. Let $M$ be given in Monge
form (\ref{eq:monge}).
As  in \S\ref{s:theory} the two cases 
(i) elliptic: $\kappa_1>0, \kappa_2>0$ and  
(ii) hyperbolic: $\kappa_1>0, \kappa_2<0$ 
are different: in the former case
the sections $f(x,y)=k$ are, locally to the origin, closed curves for small $k>0$ whereas in the latter case they are open curves extending to infinity for both signs
of $k$.

In this article our applications will concentrate on the cases where the intersection is a closed curve, 
but we give some details of other cases in this section.

$\bullet$ \ For a hyperbolic point we seek to parametrize the section by
\begin{eqnarray}
 x &=& t\;r(t,\t)\cosh\t, \ y=t\;r(t,\t)\sinh\t, \ z=t^2, \ \ \mbox{or} \nonumber \\
  x& = &t\;r(t,\t)\sinh\t, \ y=t\;r(t,\t)\cosh\t, \ z=-t^2,
\label{eq:param-hyp}
\end{eqnarray}
for a suitable function $r$. 

$\bullet$ \ For an elliptic point we seek to parametrize the section by
\begin{equation}
 x = t\;r(t,\t)\cos\t, \ \ y= t\;r(t,\t)\sin\t, \ \ z = t^2, \ t>0,\ 0\le\t< 2\pi.
\label{eq:param-ellip}
\end{equation}

$\bullet$ \ For an `elliptic cusp of Gauss' we seek to parametrize the section by
\begin{equation}
 x = (t\;r(t,\t))^2\cos\t, \ \ y= t\;r(t,\t)\sin\t, \ \ z = t^4, \ t>0,\ 0\le\t< 2\pi.
\label{eq:param-ellip-cog}
\end{equation}

We will write
\[ r(t,\t)=r_0+r_1t + r_2t^2 +\ldots,\]
where the $r_i$ are functions of $\t$ only. In the elliptic case for small values of $k>0$ the section $z=k$ is, locally to the origin, a
closed curve and the $r_i$ will
be periodic and hence functions of $\cos\t$ and $\sin\t$.

We therefore need to express the $r_i$ in terms of the Monge form coefficients $\kappa_i, \ b_i, \ c_i$
and so on. This is done by comparison of series, and enables us to approximate the surface $z=f(x,y)$,
and the sections $z=$ constant
up to any desired accuracy.
 
Consider the terms of degree $i$ in the Taylor expansion of $f$: this is a homogeneous polynomial
of degree $i$ in $x$ and $y$. 
Denote by $p_i$ the function of $\theta$ obtained from this homogeneous polynomial as follows:\\
{\em hyperbolic case}:  replace  $x$ by $\cosh\theta$  and $y$ by $\sinh\theta$ when $z>0$
and\\
 $x$ by $\sinh\theta$  and $y$ by $\cosh\theta$ when $z<0$;\\
{\em elliptic case}:
replace $x$ by $\cos\theta$  and $y$ by $\sin\theta$. 
Thus

\smallskip\noindent
$p_2=\frac{1}{2}(\kappa_1c^2+\kappa_2s^2);$
\newline
$p_3=b_0c^3+b_1c^2s+b_2cs^2+b_3s^3$;
\newline
$p_4=c_0c^4+c_1c^3s+c_2c^2s^2+c_3cs^3+c_4c^4$;
\newline
$p_5=d_0c^5+d_1c^4s+d_2c^3s^2+d_3c^2s^3+d_4cs^4+d_5s^5$;
\newline
and so on, 

\smallskip\noindent
where \\ 
{\em hyperbolic case}: $c=\cosh\t$, $s=\sinh\t$ if $z>0$ and $c=\sinh\t$, $s=\cosh\t$ if $z<0$,\\
{\em elliptic case}: $c=\cos\t$ and $s=\sin\t$.

\medskip\noindent
{\bf Hyperbolic case} \  
By substitution in the Monge form (\ref{eq:monge}), we obtain in succession the following formulas.
\newline
$r_0=\frac{1}{\sqrt{p_2}}$ when $z>0$ and $r_0=\frac{1}{\sqrt{-p_2}}$ when $z<0$ ;
\newline
$r_1=-\frac{1}{2r_0p_2}r_0^2p_3;$
\newline
$r_2=-\frac{1}{2r_0p_2}(r_0^4p_4+3r_0^2r_1p_3+r_1^2p_2)$;
\newline
$r_3=-\frac{1}{2r_0p_2}(2r_1r_2p_2+r_0^5p_5+3r_0^2r_2p_3+3r_0r_1^2p_3+4r_0^3r_1p_4)$;
\newline
$r_4=-\frac{1}{2r_0p_2}(5r_0^4r_1p_5+4r_0^3r_2p_4+2r_1r_3p_2+6r_0^2r_1^2p_4+r_2^2p_2+6r_0r_1r_2p_3$ 
\newline
\hspace*{ .9cm} $+ \ 3r_0^2r_3p_3+r_1^3p_3+r_0^6p_6)$;
\newline
$r_5= -\frac{1}{2r_0p_2}(3r_0r_2^2p_3+3r_1^2r_2p_3+6r_0r_1r_3p_3+3r_0^2r_4p_3+12r_0^2r_1r_2p_4$

\hskip .5truecm $+ \ 4r_0^3r_3p_4+4r_0r_1^3p_4+6r_0^5r_1p_6+5r_0^4r_2p_5+
10r_0^3r_1^2p_5+2r_1r_4p_2+2r_2r_3p_2)$;
etc.

\medskip

Taking 
the quadratic terms of the surface to be $x^2-\alpha^2y^2$ where $\alpha>0$, 
we have
\[r_0^2=\frac{1}{(1-\a^2)\cosh^2\t+\a^2} \ (z>0); \ \ \ \ r_0^2=\frac{1}{1-(1-\a^2)\cosh^2\t} \ (z<0).\] 

When $z>0$, the expression for $r_0^2$ is $>0$ for all $\t$ when $0 < \a \le 1$, but
for $\a>1$ we need  $\cosh^2\t< \frac{\a^2}{ \a^2-1}$, 
that is  $1\leq \cosh \theta< \frac{\a}{\sqrt{\a^2-1}}$.

When $z<0$, the expression for $r_0^2$ is $>0$ for all $\t$ when $\a\ge 1$ but if  $0<\a<1$ we need 
$1\leq \cosh\theta<\frac{1}{\sqrt{1-\a^2}}$. 

\medskip\noindent
{\bf Elliptic case} \  This is much simpler because, as above, we can expect a global parametrization of the closed
curves $f(x,y)=k$, locally to the origin, for $k>0$.  The formulae for $p_i$ are given above and those for $r_i$ are exactly the same as 
for the hyperbolic case, except that, for $z>0$ only,
\[ r_0^2=\frac{2}{\kappa_1\cos^2\theta +\kappa_2\sin^2\theta}.\]
Note that we are assuming $\kappa_1>0, \kappa_2>0$ since the surface is assumed to be
above the plane $z=0$ close to the origin.  Thus $r_0$ is always real.

When the origin is an umbilic point, then $\kappa_1=\kappa_2=\kappa$, say, and $r_0^2=\frac{2}{\kappa}$,
a constant. Indeed, we always scale the variables so that the quadratic terms of $f$ are $x^2+y^2$ and
then $r_0^2$ has the constant value 1 and we can take $r_0=1$ without loss of generality.

\medskip\noindent
{\bf Elliptic cusp of Gauss case} \ The relevant parametrization here is~(\ref{eq:param-ellip-cog}), and we apply
Proposition~\ref{prop4}. For a cusp of Gauss, after renaming axes if necessary, we have
\begin{equation}
\kappa_1\ne 0, \ \kappa_2=0,
\ b_3=0 \ \mbox{and} \left\{ \begin{array}{rl} b_2^2 < 2\kappa_1c_4 & \mbox{elliptic cusp},\\
                                                                         b_2^2>2\kappa_1c_4 & \mbox{hyperbolic cusp}. \end{array} \right. 
\label{eq:cog}
\end{equation}
Consider the elliptic case.
The expansion of the function $r$ in Proposition~\ref{prop4}  is obtained by a `weighted' version of the
method used in the 
previous cases. This time let $p_i, \ i\ge 4$, be the result of substituting $x=\cos\theta, \ y=\sin\theta$
in the terms of $f$ of weighted degree $i$, where $x$ has weight 2 and $y$ has weight 1. Thus writing $c=\cos\theta,
\ s=\sin\theta$ we have\\
$p_4=\frac{1}{2}\kappa_1c^2+b_2cs^2+c_4s^4 = q(c,s) $ (see Proposition~\ref{prop4})\\
$p_5= b_1c^2s+c_3cs^3+d_5s^5$\\
and so on. Writing $r=r_0+r_1t+r_2t^2+r_3t^3+\ldots$ as before, where each $r_i$ is a function of $\theta$ only, we find on
substitution\\
$r_0^4=1/p_4$,\\
$4p_4r_1+r_0^2p_5=0$, from which we solve for $r_1$, \\
$4p_4r_0r_2+6p_4r_1^2+5p_5r_0^2r_1+p_6r_0^4=0$, from which we solve for $r_2$, \\
and so on.  The assumption of an {\em elliptic} cusp of Gauss guarantees that $p_4$ is nonzero for
all values of $\theta$.

\medskip

The same ideas can be used in principle to parametrize the sections near a hyperbolic cusp of Gauss. It is
not difficult to show that in this case (where $b_1^2>2\kappa_1c_4$) the expression $p_4$ above
vanishes for exactly four values of $\theta$ in the range $0< \theta < 2\pi$. In increasing order these
are of the form $0<\theta_1<\theta_2<\theta_3=2\pi-\theta_2<\theta_4=2\pi-\theta_1<2\pi$. In the ranges from
$\theta_1$ to $\theta_2$ and from $\theta_3$ to $\theta_4$ the expression $p_4$ is negative and in the
other two ranges it is positive. Then we can use \\
$x=t^2r^2\cos\theta, \ y=tr\sin\theta, \ z=-t^4$ in the first two ranges and \\
$x=t^2r^2\cos\theta, \ y=tr\sin\theta, \ z=t^4$ in the other two ranges, \\
determining $r$ exactly as before, except that $r_0^4=-1/p_4$ when $p_4<0$.  Needless to say this
case is much more delicate than the elliptic case since the branches are not closed and the relationship
between the value of $t$ and the closeness of fit obtained by taking say four terms of the
approximation to $r$ will be difficult to determine.  An example is shown in Figure~\ref{fig:approx}.

\section{Applications to pre-symmetry sets and symmetry sets}\label{s:appl}

Let $\ga$ be a smooth {\em parametrized} plane curve. We seek first the `pre-symmetry set' which is the set of parameter pairs
$(u_1,u_2)$ such that there exists a circle tangent to $\ga$ at $\ga(u_1), \ga(u_2)$. The pre-symmetry
set is symmetric about the `diagonal' $u_1=u_2$. A convenient way to find these pairs
is to look for solutions of the equation given by a scalar product 
\begin{equation}
(\ga_1 - \ga_2)\cdot(\T_1-\T_2)=0,
\label{eq:pre-ss}
\end{equation}
where $\ga_i=\ga(u_i)$ and $\T_i$ is the unit oriented tangent to $\ga$ at $\ga_i$.  These solutions are precisely the pairs required, together
with any pairs where $\T_1=\T_2$: parallel oriented tangent pairs. For a closed convex curve the only such pairs
arise from $u_1=u_2$ and these are easily identified as the diagonal. 
When a closed curve has inflexions other pairs with $\T_1=\T_2$ will be present. 
We need to beware of these when interpreting the pre-symmetry set.
The pre-symmetry set contains much information about the symmetry set. Crossings of the diagonal correspond to endpoints
of the symmetry set and horizontal or vertical tangents correspond to cusps. 

Secondly we want to plot the centres of the bitangent circles identified above. There are many ways to do this but here
is a convenient one, couched in the language of complex numbers. Let $c$ be the centre of the circle. Then, for some angle
$\theta$,
\[ (\ga_2-c)=e^{i\t}(\ga_1-c), \ \ \T_2=e^{i\t}\T_1.\]
Then, as complex numbers, we have by dividing these equations $c=(\ga_2\T_1-\ga_1\T_2)/(\T_1-\T_2)$.
Note that if $\T_1=\T_2$, or, in practice, if these are sufficiently close together, and $\ga_1\ne \ga_2$, then $c$ is very far away
and it will not appear on the diagram of the symmetry set.

The above has been implemented in the Liverpool Surfaces Modelling Package (LSMP, \cite{lsmp}), also known
as SingSurf. Note that it is necessary for $\ga$ to be {\em parametrized}, which is the motivation for the
theoretical results given in \S\ref{s:theory}.
We give below examples of umbilic points and elliptic cusps of Gauss (the latter obtained
for the present by a different method).  We shall give examples of hyperbolic points
elsewhere. See~\cite{diatta-giblin2004,scale-space2005} for other methods and theoretical results
related to the examples below.

\subsection{Umbilic points}\label{s:umbilic}
Elliptic points for which the principal curvatures $\kappa_i$ as in (\ref{eq:monge}) are unequal make rather uninteresting
examples since the intersection $f(x,y)=k$ is a curve which is very nearly an ellipse, having exactly four vertices.
The symmetry set has two smooth branches and the medial axis has one smooth branch, as for an ellipse. We concentrate
here on the case of an {\em umbilic} point, where $\kappa_1=\kappa_2$, and without loss of
generality we can take these both equal to 1 by scaling the surface. Thus the equation of the surface has
the form 
\[ z=x^2+y^2+b_0x^3+b_1x^2y+b_2xy^2+b_3y^3+ \ \mbox{h.o.t.}\]

\begin{figure}
  \begin{center}
  \leavevmode
\epsfysize=1.9in
\epsffile{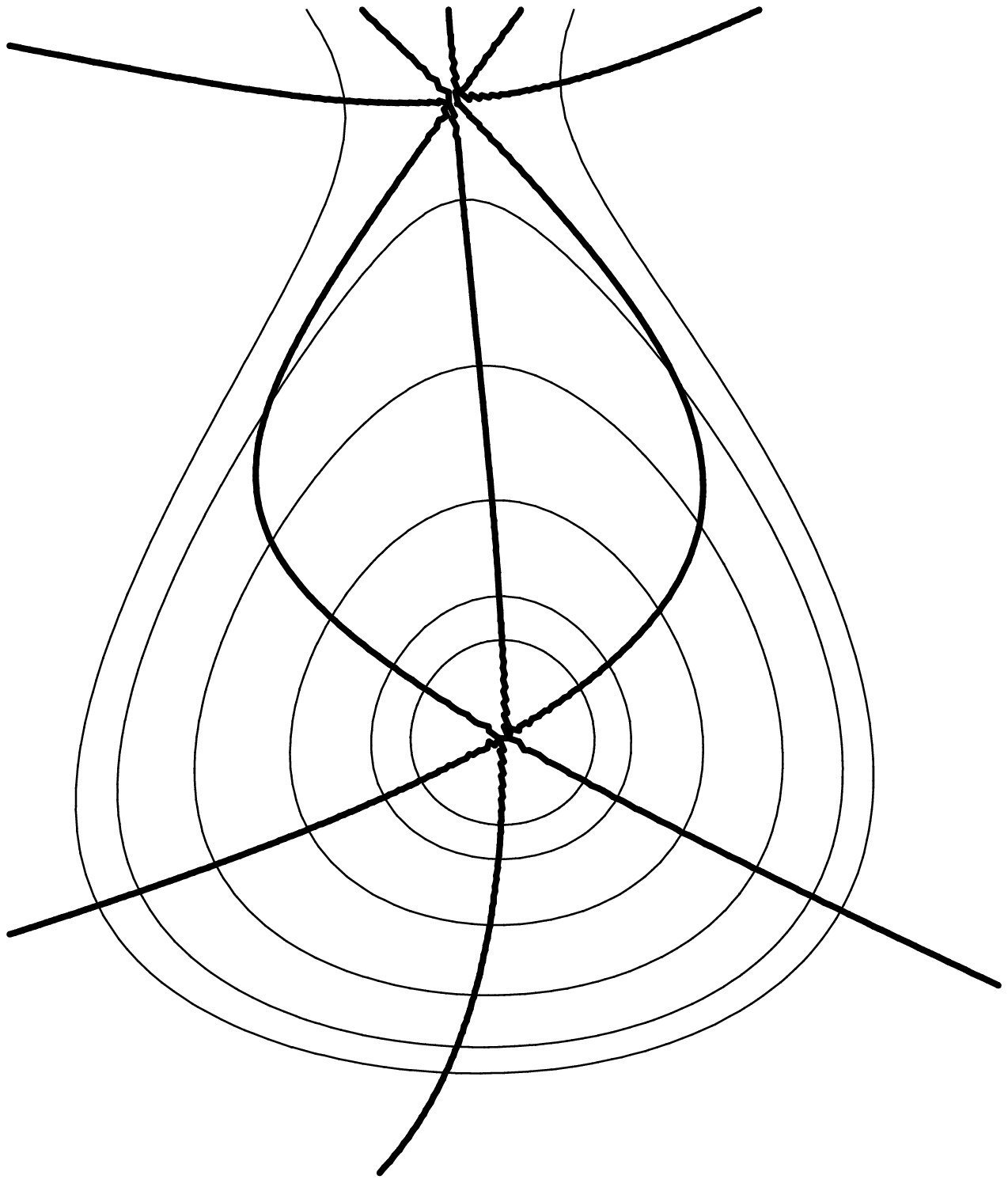}
\epsfysize=2.4in
\epsffile{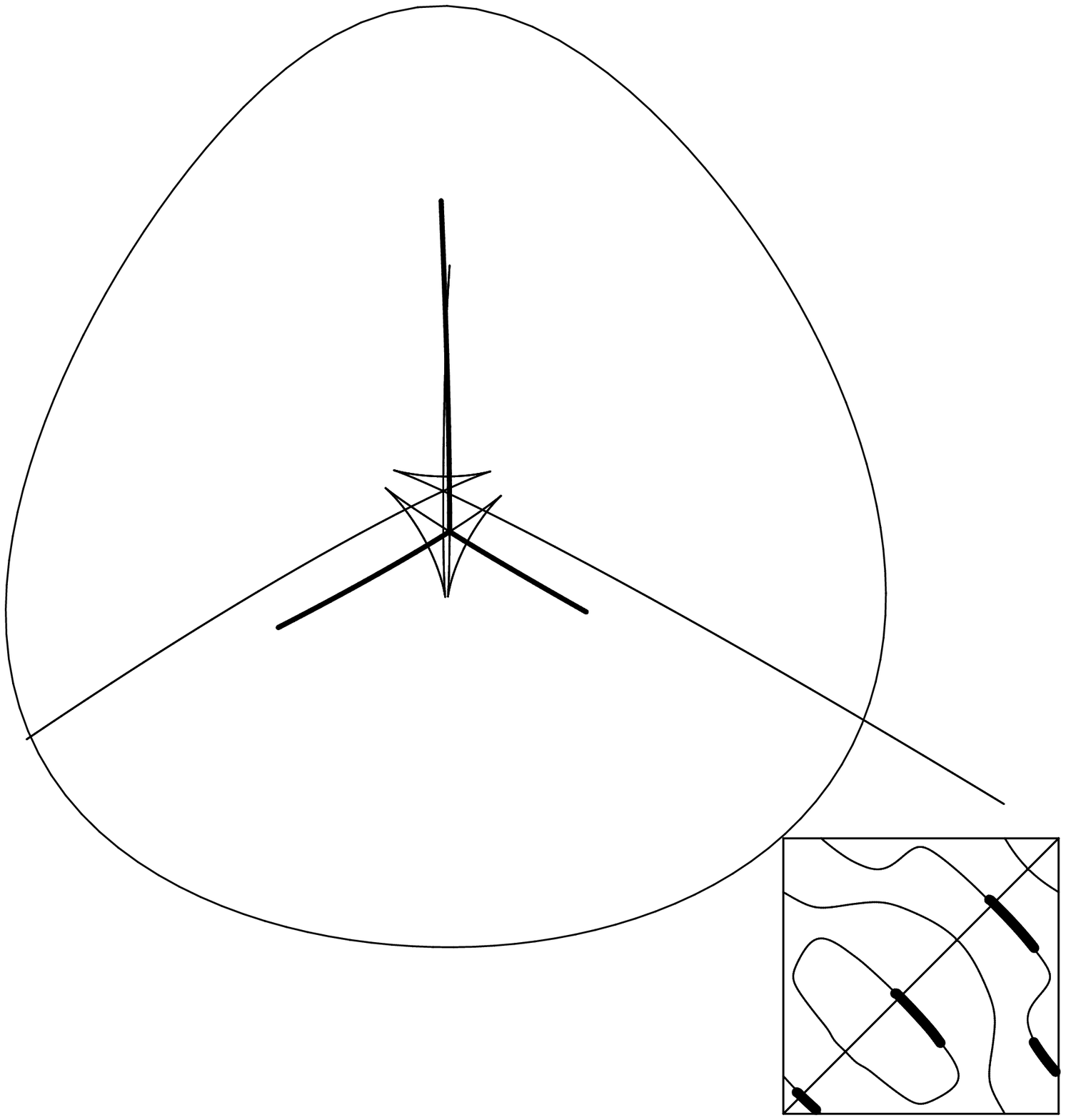}
\end{center}
\vspace*{-0.5in}
\begin{center}
\epsfysize=1.9in
\epsffile{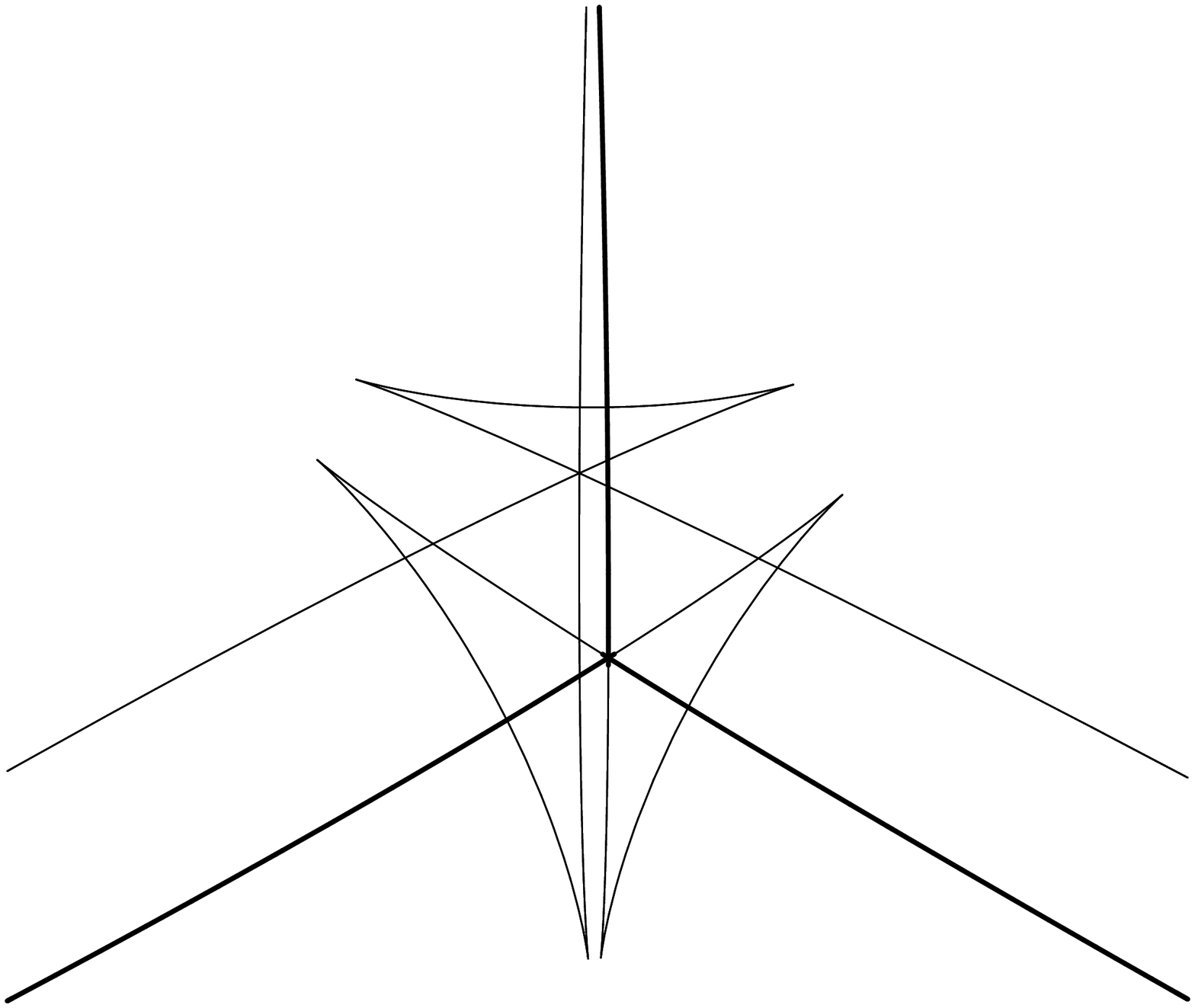}
\epsfysize=2.4in
\epsffile{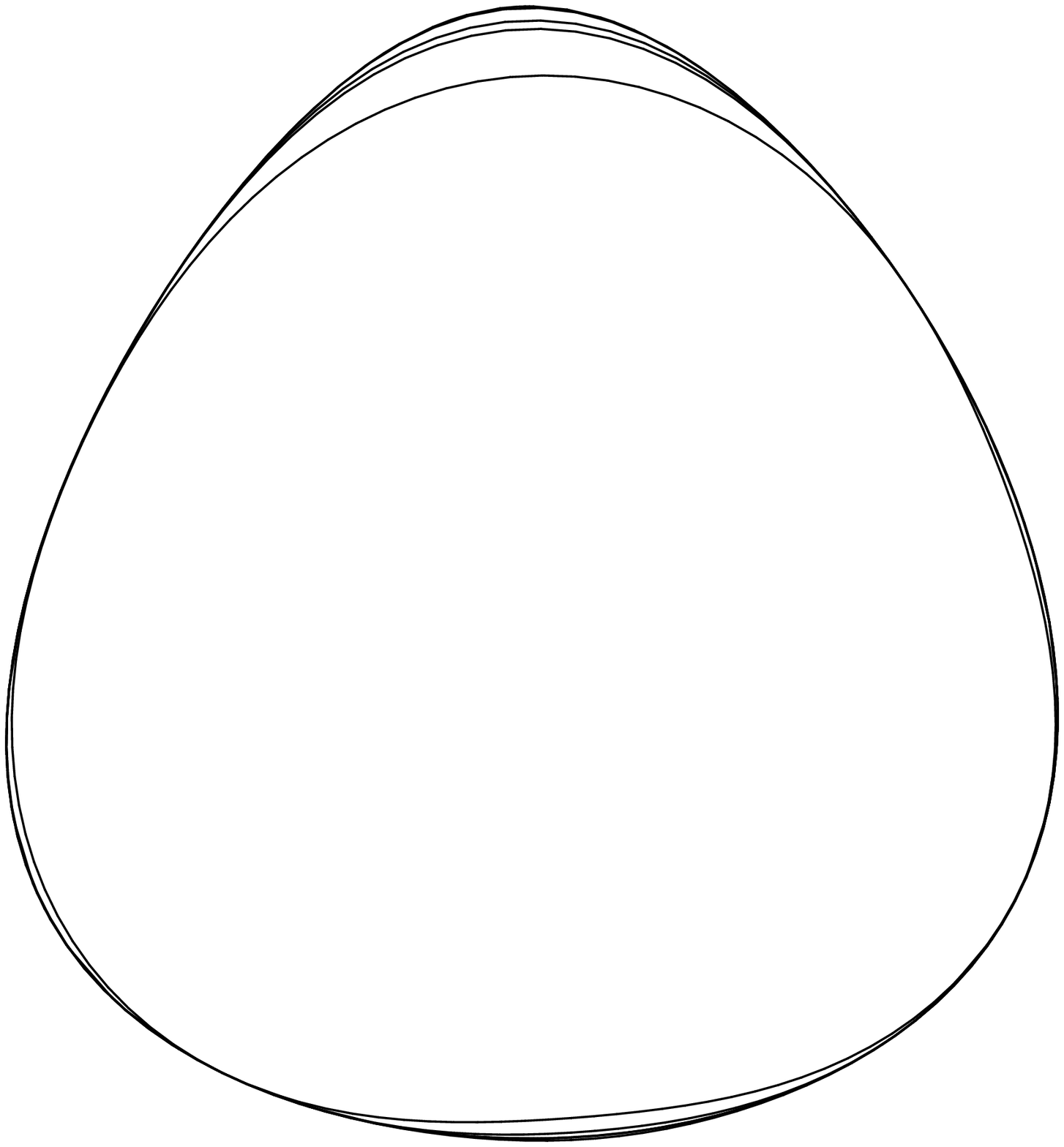}
 \end{center}
\vspace*{-0.5in}
\caption{The umbilic case. Top left: several level sets $f(x,y)=x^2+y^2+x^3+2x^2y+xy^2-3y^3=k$ for $k$ from 0.001 to 0.018.
The curve becomes more circular as $k$ decreases.
The thick curves are the locus of vertices, obtained as in \protect{\cite{diatta-giblin2004}}. Notice that if $k$ is too big
then the level set ceases to be
closed. Top right:  The symmetry set and (boxed insert) the pre-symmetry set of the
curve $f=k$ for a small enough value of $k$ (0.01 in fact) that the structure has stabilised. The heavily drawn part of the
pre-symmetry set corresponds to the Y-shaped medial axis. The curve $f=k$ was parametrized by
the method of \S\ref{s:practice}, using 10 terms.  Bottom left: A closeup of the central part of the symmetry set; note the
six cusps and two triple crossings. The medial axis is
drawn heavily. Bottom right: The approximations  to $f(x,y)=0.01$ given by taking 1, 3, 4 and 10 terms as
in \S\ref{s:practice}: the outermost curve is for 10 terms and is essentially indistinguishable from the exact solution.}
\label{fig:umbilic}
 \end{figure}

 By a technique explained in
\cite{diatta-giblin2004} we can find the locus of vertices of the family of curves $f(x,y)=k$ for $k>0$, and in
\cite{scale-space2005} there are results on triple intersections and cusps of the symmetry set. We shall not recall these
in detail here but will point out how they are verified by the example.
Clearly any rotation about the $z$-axis will leave the equation of the surface
in the same form as above. If we rotate to make $b_0=b_2$
then it can be shown that the three branches of the vertex locus make angles with the $x$-axis which are
multiples of $60^\circ$, as in Figure~\ref{fig:umbilic}, top left. We pause here to establish the existence of such a rotation.

\begin{lemma}\label{rotatecoordinates} Let $f$ be given by  $f(x,y)=x^2+y^2+b_0x^3+b_1x^2y+b_2xy^2+b_3y^3+$ h.o.t.  
 Then we can rotate the coordinates in the $x,y$ plane so that in the new coordinates $u,v$, the function $f$ takes the form
 $f=u^2+v^2+B_0u^3+B_1u^2v+B_0uv^2+B_3v^3+$ h.o.t.\end{lemma}
\proof

We rotate by an angle $\phi$. Expressing $f$ in terms of the new 
coordinates $(u,v)$, amounts to replacing $x$ by $u\cos\phi+v\sin\phi$ and $y$ by $-u\sin\phi+v\cos\phi$ in the expression for $f(x,y)$. 
The new expression for $f$ is of the form \\
$f=u^2+v^2+B_0u^3+B_1u^2v+B_2uv^2+B_3v^3+$ h.o.t.,   
where\\
$B_0= b_0\cos^3\phi-b_1\cos^2\phi\sin\phi+b_2\cos\phi\sin^2\phi-b_3\sin^3\phi,$\\
$B_2= 3b_0\cos\phi\sin^2\phi-b_1\sin^3\phi+2b_1\cos^2\phi\sin\phi-2b_2\cos\phi\sin^2\phi$\\
\hspace*{0.5in} $+b_2\cos^3\phi-3b_3\cos^2\phi\sin\phi.$

We need to show the existence of an angle $\phi_0$ for which $B_0=B_2$.
  Substitute $U=\tan(\frac{\phi}{2})$ and write $p=b_2-b_0, \ q=b_3-b_1$. Then the equation $B_0=B_2$ reads
\[ p U^6 + 6q u^5-15p U^4-20q U^3+15p u^2+6q U - p = 0.\]
Note that the equation has pairs of roots of the form $U,-1/U$, corresponding to solutions
$\phi, \phi+\pi$. Of course, if $p=0$ we can take $\phi=0$. Otherwise the discriminant of this degree 6 equation
is a positive constant times $(p^2+q^2)^5$ and hence $>0$. For a degree six equation this means that there
are either two or six real solutions for $U$. \qed

\smallskip\noindent
{\bf Remark}
We do not know a geometrical interpretation of the distinction between two and six solutions here.
Note that these correspond to respectively one and three pairs of solutions $\phi, \phi+\pi$. There are
a number of situations where through an umbilic pass one or three geometrically defined curves,
for example ridges and sub-parabolic curves; see \cite{b-g-t96b}. The present case appears
to be different from these.

\smallskip

In Figure~\ref{fig:umbilic} we show an example to illustrate the above methods. It confirms the results
of \cite{diatta-giblin2004} and \cite{scale-space2005}, for small enough $k>0$,  namely: (a) there are six vertices on the level set,
resulting in three branches, of which one connects a maximum to maximum of curvature, one a minimum to minimum and
one a maximum to minimum,
 (b) there are two triple crossings on the symmetry set, (c) there are six cusps on the symmetry set. Note that full detail
is given on the symmetry set since the polar approximation  is smooth and extremely close to the level set $f=k$.

\begin{figure}[h!]
  \begin{center}
  \leavevmode
\epsfysize=2in
\epsffile{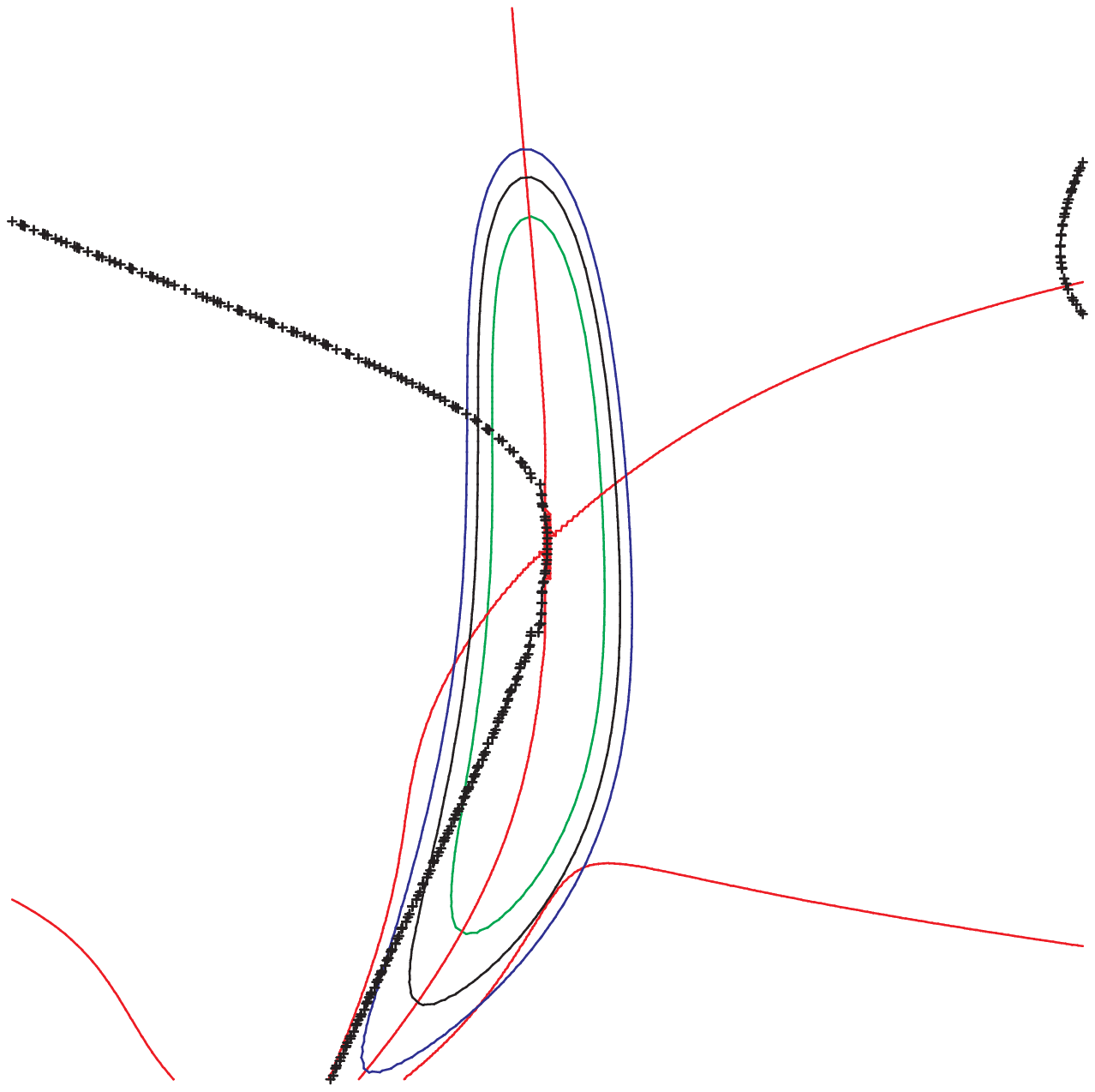}
\hspace*{0.5in}
\epsfysize=2in
\epsffile{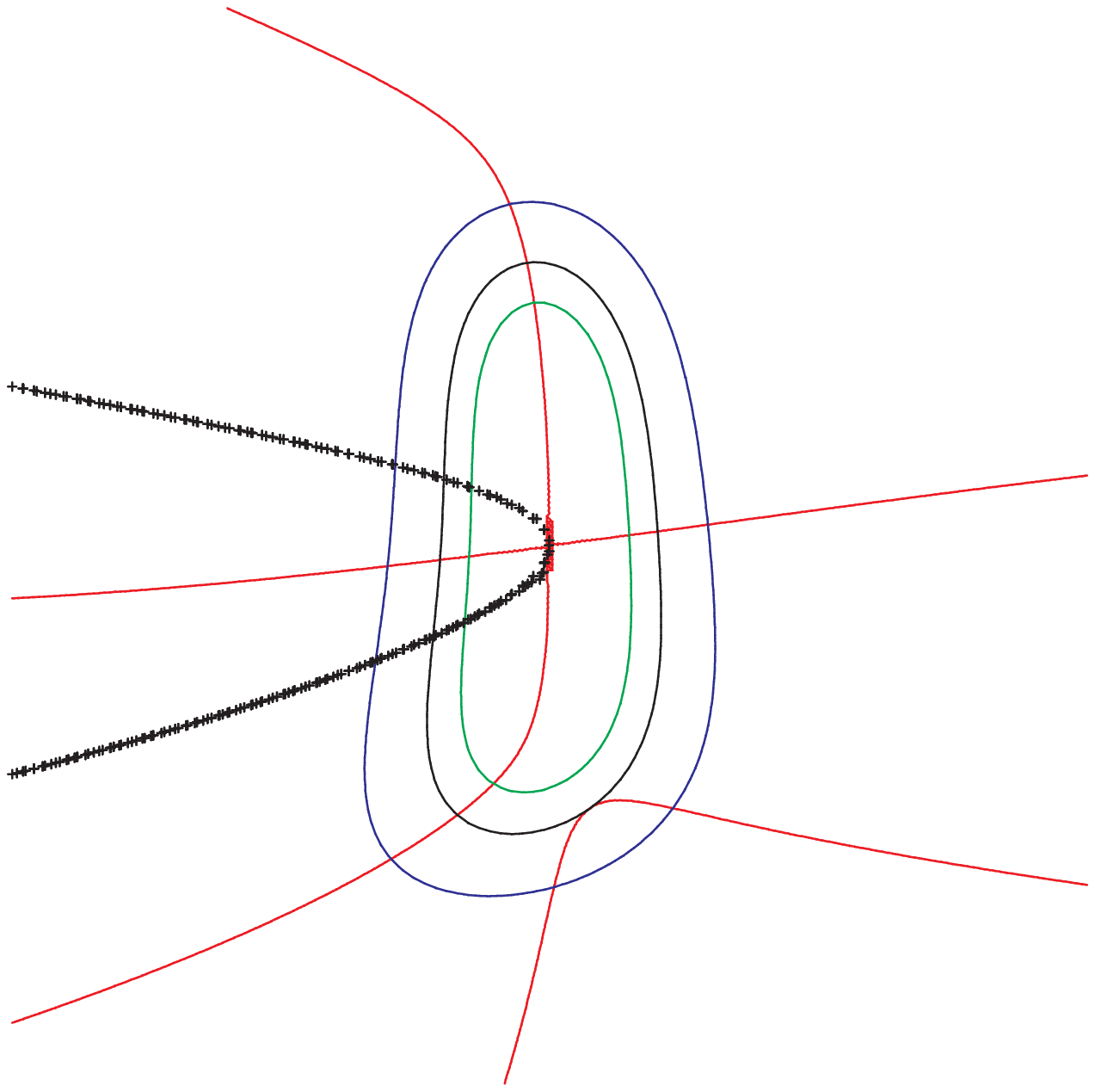}
 \end{center}
\caption{Elliptic cusps of Gauss, examples 1 (left) and 2 (right). The closed curves are the level sets $f(x,y)=k$, the other
thin curves are the loci of vertices of  level sets and the thick curves are the loci of inflexions. Notice that for small
enough $k$ there are four vertices (as predicted by \cite{diatta-giblin2004}) but as $k$ increases the level set is first
tangent to a branch of the locus of vertices not local to the origin and then intersects this locus in two points, giving six vertices.}
\label{fig:cog1}
  \end{figure}

\subsection{Cusps of Gauss}
\label{s:cog}
We consider here elliptic cusps of Gauss, as in (\ref{eq:cog}), and assume $\kappa_1>0$ so that the level set $f=k$ is, locally
to the origin, a simple closed curve for small $k>0$.  According to the general results of \cite{diatta-giblin2004} there are four
vertices and two inflexions on this closed curve. (The situation for hyperbolic cusps of Gauss---and indeed for
hyperbolic points---is much more complicated,
with the pattern of vertices and inflexions giving rise to several cases, which are detailed in \cite{diatta-giblin2004}.)
Figure~\ref{fig:cog1} shows several level sets together with the loci of vertices and inflexions, for two elliptic
cusps of Gauss. \\
Example 1 is $f(x,y) = x^2+2x^2y+xy^2+2x^2y^2-xy^3+y^4$, \\
Example 2 is $f(x,y) = x^2+2x^2y+xy^2+2x^2y^2-xy^3+6y^4$.

\begin{figure}
  \begin{center}
  \leavevmode
\epsfysize=1.3in
\epsffile{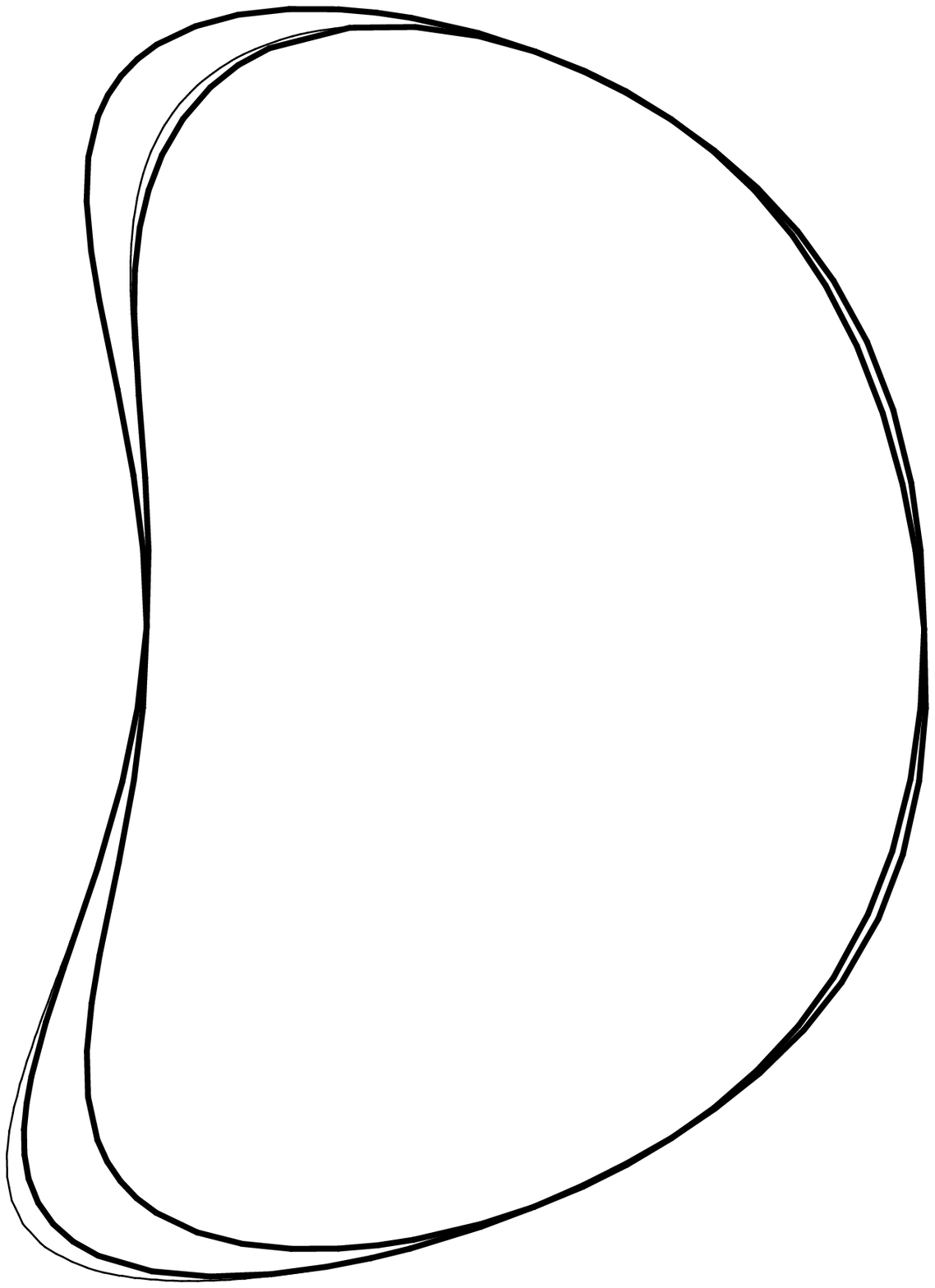}
\epsfysize=1.5in
\hspace*{0.7in}
\epsffile{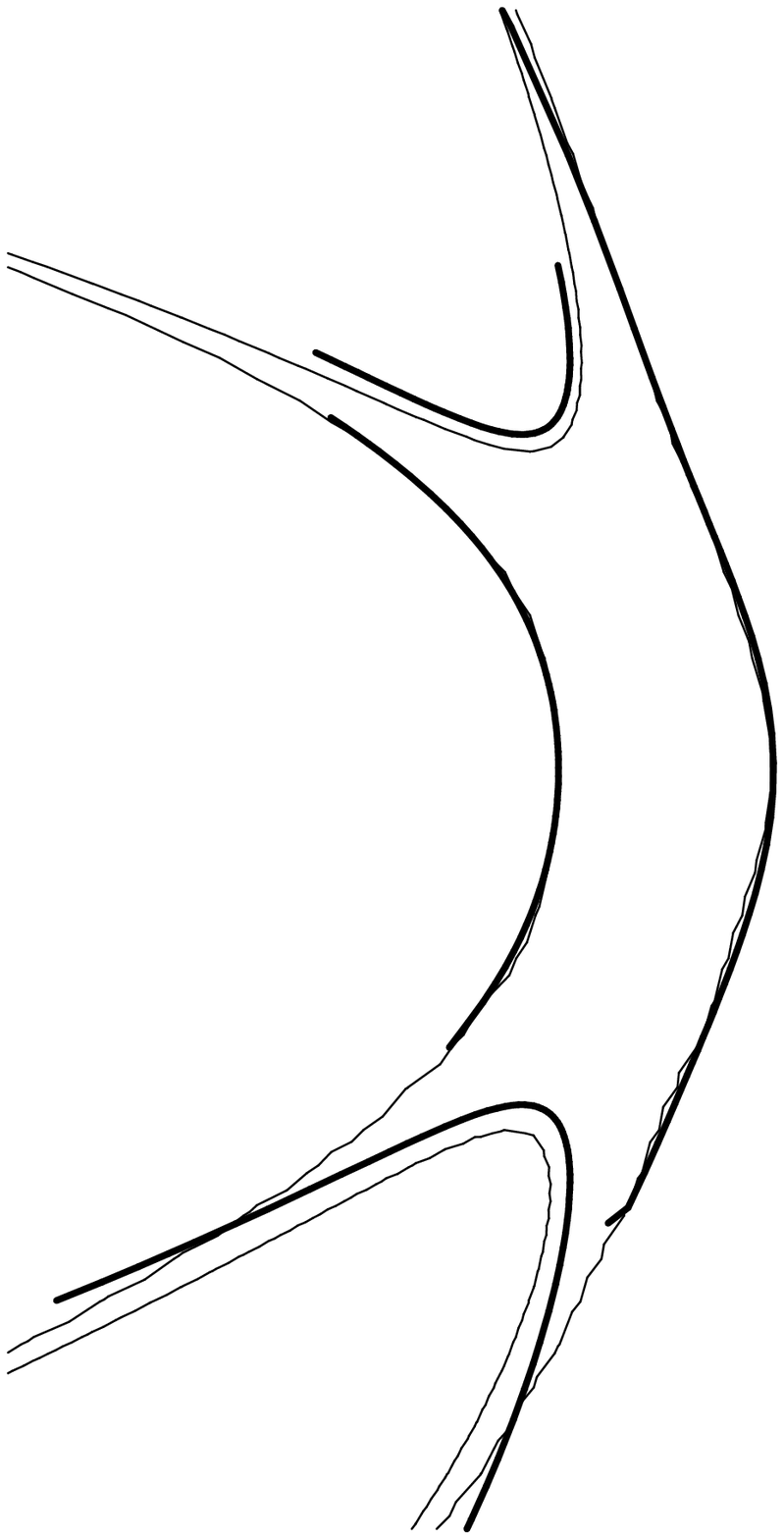}
\end{center}
\caption{Left: the thin curve is a level set $f(x,y)=k$ for an elliptic cusp of Gauss
as in Example~1 in the text, and the thick lines are obtained
by taking $r=r_0$ and $r=r_0+r_1t$ in the approximations described in \S\ref{s:practice}. The approximation
obtained by taking three terms $r=r_0+r_1t+r_2t^2$ is visually indistinguishable from the true level set. Right: the thin curves are
two level sets $f(x,y)=k$ for values of $k$ with opposite signs in the case of a hyperbolic cusp of Gauss
$f(x,y)= x^2+3xy^2+3xy^3+y^4$. The thick
curve is obtained with the first term $r=r_0$ of the approximation of \S\ref{s:practice}. Clearly there are more problems
here with capturing all the local features of the level sets. {\em Note:} \ in these figures the horizontal scale
has been exaggerated compared to the vertical scale, as an aid to clarity.} 
\label{fig:approx}
 \end{figure}

Figure~\ref{fig:approx}, left, shows the curve of Example~1 together with approximations obtained
as in \S\ref{s:practice}. In this figure the horizontal scale is exaggerated to improve clarity (otherwise the approximations
and the actual curve are too hard to distinguish!). We also show a hyperbolic cusp of Gauss example in Figure~\ref{fig:approx}, right, though we
do not go on to examine the symmetry set here.

\begin{figure}
  \begin{center}
  \leavevmode
\epsfxsize=4in
\epsffile{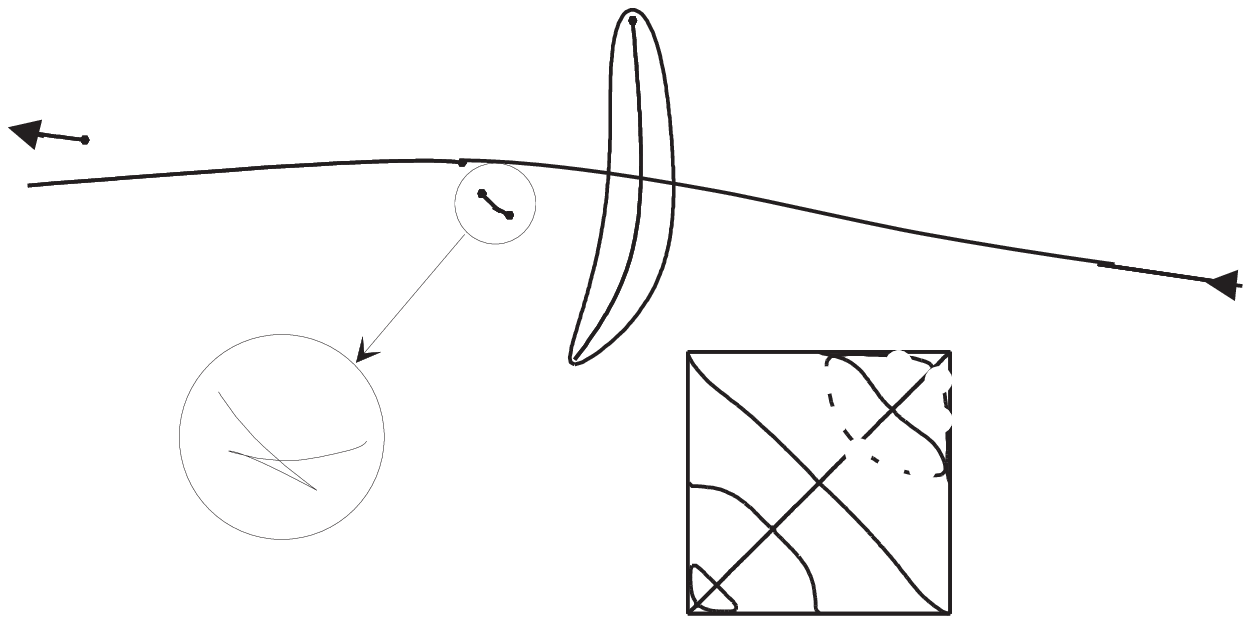}
\end{center}
\vspace*{-2.5in}
\begin{center}
\epsfxsize=3.5in
\epsffile{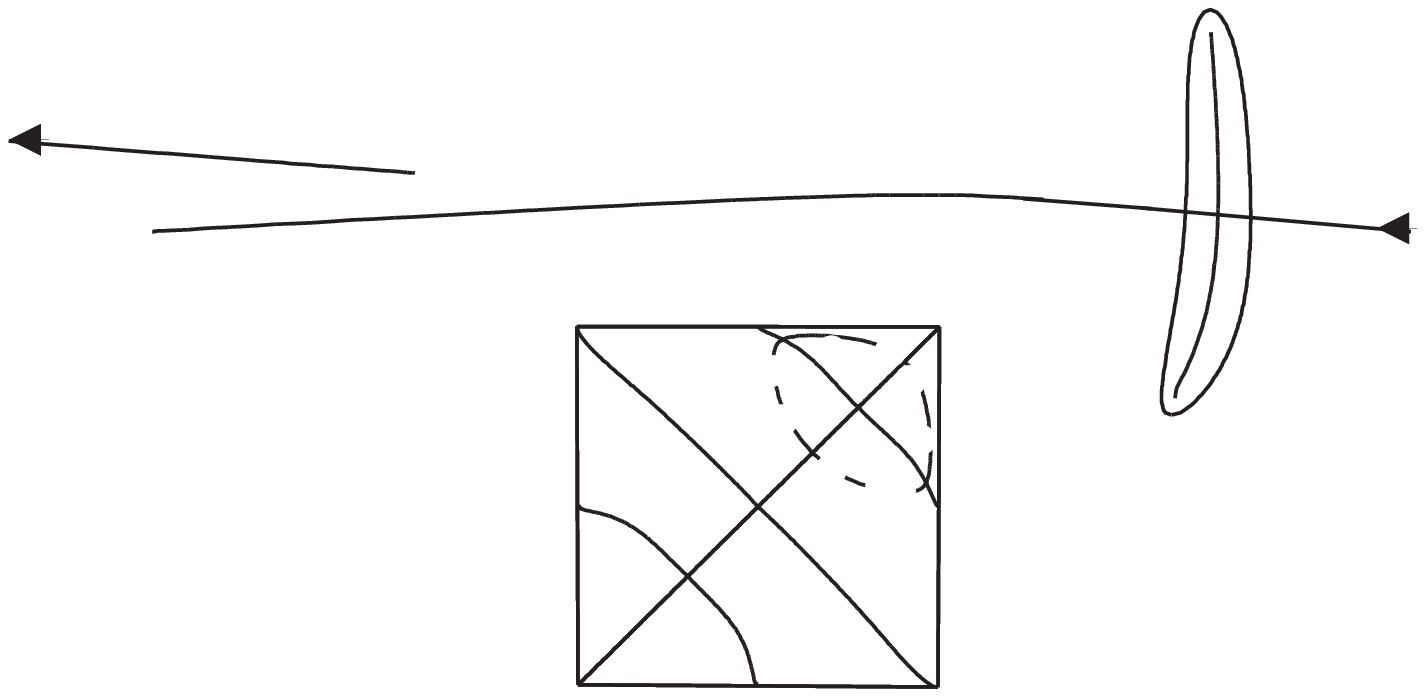}
 \end{center}
\vspace*{-3in}
\caption{Example 1 as in Figure~\ref{fig:cog1}, with annotated symmetry set and pre-symmetry set. The dashed section
of the pre-symmetry set (boxed inserts) is due to the parallel tangents arising from the inflexions on $f(x,y)=k$; see \S\ref{s:appl}. It is
ignored when plotting the symmetry set (main diagram). The arrows indicate that a branch of the symmetry set goes to infinity when
the bitangent circle becomes a bitangent line. The other endpoints are actual. Top: $k$ is large enough for the level set to
have six vertices. The additional branch of the pre-symmetry set is the loop towards the bottom left of the
box and the corresponding piece of symmetry set is enlarged in the circle. 
As $k$ decreases the number of vertices reaches its stable value of 4, as in the bottom figure. 
The medial axis for both values of $k$ consists of the branch interior to the curve and the part of the other branch going
to infinity to the left.}
\label{fig:cog2}
  \end{figure}

\begin{figure}
  \begin{center}
  \leavevmode
\epsfysize=4in
\epsffile{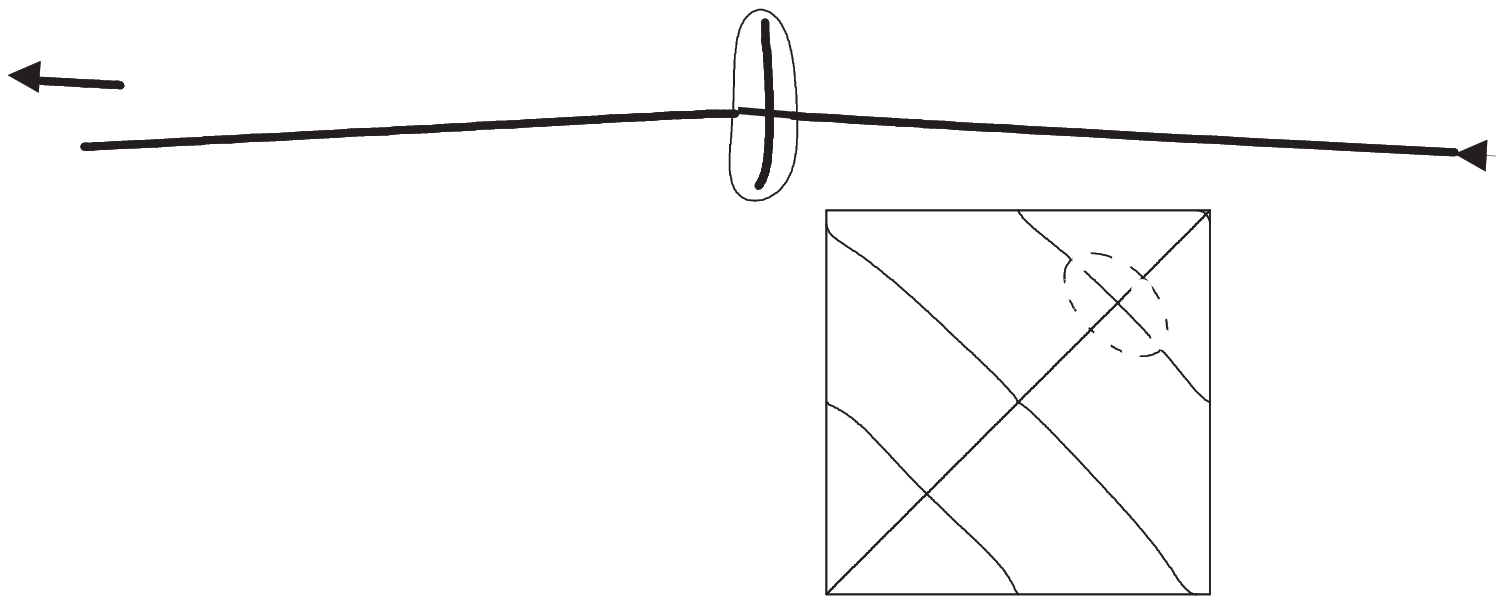}
 \end{center}
\vspace*{-1.3in}
\caption{Example 2 of Figure~\ref{fig:cog1}. Here $k$ is taken small enough to make the level set have its stable
number of vertices, namely 4. The symmetry set has very much the same structure as Example 1, that is two
simple branches one of which extends to infinity.}
\label{fig:cog3}
  \end{figure}
Theoretical results on the symmetry sets of these curves appear to be much more difficult to obtain than those
for elliptic, hyperbolic and ordinary parabolic points (\cite{diatta-giblin2004,scale-space2005}). Thus
evidence gathered from examples is all the more valuable. We show in Figure~\ref{fig:cog2} the symmetry sets
and pre-symmetry sets for the two examples of Figure~\ref{fig:cog1}. In both cases the stable situation as
$k\to 0$ is two smooth branches, one of which extends to infinity. Thus the structure of the symmetry set
for an elliptic cusp of Gauss appears to be very much simpler than for an umbilic point. In particular there
are no cusps and no triple crossings. Furthermore, the only distinction between the symmetry sets of sections of
a surface close to the tangent plane at an ordinary elliptic point and at an elliptic cusp of Gauss is that in the
latter case one of the branches of the symmetry set extends to infinity. On the other hand, the distinction between
an ordinary elliptic point and an umbilic is very striking: in the latter case there are six cusps and two triple
crossings as in Figure~\ref{fig:umbilic}.

\section{Conclusion}
We have given a simple method of parametrizing closed level set curves $f(x,y)=k$ to
any desired degree of accuracy and used it to compute symmetry sets
and medial axes
of sections of a surface $z=f(x,y)$ in Monge form close to the tangent plane at an umbilic and an elliptic cusp of Gauss.
 The method applies also to level sets $f(x,y)=k$ which are not closed; here the patterns of vertices and inflexions
on  the plane sections become much more complicated \cite{diatta-giblin2004} and there are many different
cases. Furthermore there are technical problems in ensuring that the approximations capture
all the necessary local information about the level sets.
We shall pursue these cases elsewhere. 

There are many related questions; for example, when a point on
a surface moves from the hyperbolic region to the parabolic curve, how does the family of symmetry sets
of the parallel plane sections behave? This is tantamount to considering a one-parameter family of {\em surfaces},
that is the curves will now belong to a {\em two}-parameter family. This again is the subject of further work.

\end{document}